\documentclass{amsart}
\usepackage{amssymb}
\usepackage{color}
\definecolor{darkgreen}{cmyk}{1,0,1,.2}
\definecolor{m}{rgb}{1,0.1,1}
\definecolor{green}{cmyk}{1,0,1,0}

\definecolor{test}{rgb}{1,0,0}
\definecolor{cmyk}{cmyk}{0,1,1,0}

\newtheorem{Equation}{}[section]

\newtheorem{theorem}[Equation]{Theorem}

\newtheorem{corollary}[Equation]{Corollary}

\def\ch{\operatorname{ch}}

\def\cohom{\operatorname{H}}

\def\Ind{\operatorname{Ind}}

\def\C{\mathbb C}

\def\R{\mathbb R}

\def\Z{\mathbb Z}
\def\T{\mathbb T}

\def\Q{\mathbb Q}

\def\what{\widehat}
\def\wtit{\widetilde}

\def\wL{\widetilde{L}}

\def\dd{\displaystyle}
\def\pa{\partial}

\begin{document}



\title[An interesting example for spectral invariants]
{An interesting example for spectral invariants\\  \today}


\author[M-T. Benameur]{Moulay-Tahar Benameur}
\address{I3M, UMR 5149 du CNRS, Universit\'{e} de Montpellier 2,  34095 Montpellier, France}
\email{moulay.benameur@math.univ-montp2.fr}

\author[J.  L.  Heitsch]{James L.  Heitsch}
\address{Mathematics, Statistics, and Computer Science, University of Illinois at Chicago} 
\address{ Mathematics, Northwestern University}
\email{heitsch@math.uic.edu}

\author[C.  Wahl]{Charlotte Wahl}
\address{Leibniz-Archiv, Waterloostr. 8, 30169 Hannover, Germany}
\email{wahlcharlotte@gmail.com}



\begin{abstract} In \cite{HL1999}, the heat operator of a Bismut superconnection for a family of generalized Dirac operators is defined along the leaves of a foliation with Hausdorff groupoid.  The Novikov-Shubin invariants of the Dirac operators were assumed greater than three times the codimension of the foliation.   It was then shown that the associated heat operator converges to the Chern character of the index bundle of the operator.  In \cite{BH2008},  this result was improved by reducing the requirement on the Novikov-Shubin invariants to one half of the codimension.  In this paper, we construct examples which show that this is the best possible result.
\end{abstract}

\maketitle


\section{Introduction}

In \cite{HL1999},  Connor Lazarov and the second author gave general conditions which guarantee that the Bismut superconnection formalism extends in full generality to families of non-compact manifolds.   This allowed us to prove a families index theorem for generalized Dirac operators defined along such families. Special cases 
give the Atiyah-Singer families index theorem, the Atiyah $L^2$ index theorem \cite{Atiyah:1976}, and the Connes foliation cohomology index theorem \cite{Connes:1986}.  In addition,  we got a new theorem for fiber bundles which is a combination of the first two above, namely an $L^2$ families index theorem.   Note that  Connes and Skandalis have proven a families index theorem for families of elliptic operators defined along the leaves of a foliation. See \cite{Connes-Skandalis:1984} and \cite{Connes:1987}.   One of the conditions we required is that the Novikov-Shubin invariants, \cite{NS86a, NS86b, HL1994}, of the Dirac operators were greater than three times the codimension of the foliation.   

In \cite{BH2004, BH2008}, the first two authors extended these results to generalized Dirac operators $D$ along the leaves of a Riemannian foliation. We assumed that the projection onto the kernel of $D$ is transversely smooth, and that the spectral projections of $D^{2}$ for the intervals $(0,\epsilon)$ are transversely smooth, for $\epsilon$ sufficiently
small.    We defined the Connes-Chern character of $D$ in the Haefliger cohomology of the foliation.  We then showed that the pairing of this Connes-Chern character with a given Haefliger  $2k-$current is the same as that of the Haefliger Chern character of the index bundle of $D$, whenever the Novikov-Shubin invariants of $D$ are greater than $k$.    In particular, if the Novikov-Shubin invariants are greater than half the codimension of $F$, they are always the same.  We conjectured that this theorem is still true provided only that the Novikov-Shubin invariants are positive.  In this paper, we show that this is false, and that the result in \cite{BH2008} is the best possible, under the given hypotheses.  It is an interesting question as to what additional conditions need to be imposed for the conjecture to be true.

\section{A bit of background}

We assume that the reader is familiar with the papers \cite{BH2004, BH2008}, in particular with the concepts of those papers, including: Haefliger cohomology; generalized Dirac operators; the various Connes-Chern characters used there; the Novikov-Shubin invariants; and transverse smoothness of leafwise operators.  Suppose that $\what{D}$ is a generalized Dirac operator defined along the leaves of a Riemannian foliation $F$ of codimension $n$ on a smooth manifold $M$.   Denote by $D$ the induced leafwise operator for the foliation $F_s$ (whose leaves are the inverse image of points on $M$ under the source map from the groupoid to $M$) on the holonomy (or homotopy) groupoid of $F$.  Denote by $P_0$ the leafwise projection onto the kernel of $D$, and by $P_{\epsilon}$ the leafwise spectral projection for $D^{2}$ for the interval $(0,\epsilon)$.  

\begin{theorem}[Theorem 4.2 of \cite{BH2008}]  \label{main} 
Assume that $P_0$,  and (for $\epsilon$ sufficiently small) $P_{\epsilon}$
are transversely smooth.  For a fixed integer $k$ with $0\leq k\leq n/2$, assume that the Novikov-Shubin invariants of $D$ are greater than $k$.  
Then the $k^{th}$ component of the Chern character of the K-theory index 
of $D$ equals the $k^{th}$ component of the Chern character of the index bundle of $D$, that is, in the Haefliger cohomology $\cohom_c^{2k}(M/F)$ of $F$,
$$
\ch^k_{a}(\Ind_{a}(D)) =  \ch^k_{a}([P_{0}]).
$$
\end{theorem}  

\begin{corollary}[Theorem 4.1 of \cite{BH2008}]\label{limit}
If the Novikov-Shubin invariants of $D$ are greater than $n/2$,   then 
$$
\ch_{a}(\Ind_{a}(D)) =  \ch_{a}([P_{0}]).
$$
\end{corollary}

\section{The examples}

We now show that Corollary \ref{limit} is the best possible.

Consider the product of two $n$-tori  $\T^n \times \T^n = \R^n \times \R^n/ \Z^n \times \Z^n$, and a foliation $F$ on it which is the product of constant irrational slope foliations on the individual $\T \times \T$.  We assume that the slopes are not rationally related.  We give $\T^n \times \T^n$ the metric it inherits from the usual metric on  $\R^n \times \R^n$.  The leaves of $F$ are then just copies of $\R^n$ with the usual metric.   $F$ is a Riemannian foliation, and its holonomy and homotopy groupoids are just the product  $\T^n \times \T^n \times \R^n$, with the foliation $F_s$ given by the $\R^n$ factors.  The Haefliger cohomology of $F$ is highly non-trivial, but still computable, and $\T^n \times \T^n$ has an interesting leafwise flat bundle over it with which to pair the leafwise Dirac operator. 

Denote by $\what{E}_n \to \T^n \times \T^n$ the Hermitian bundle which is given as follows.   Let $(\xi,\what{\xi}) \in \Z^n \times \Z^n = \pi_1(\T^n \times \T^n)$, and $(x, w, z) \in  \R^n \times \R^n \times \C$,  and define  
$$
(\xi,\what{\xi}) \cdot (x, w, z) =  ( x + \xi, w + \what{\xi},  e^{2\pi i \langle \xi, w \rangle} z),
$$
where 
$$
\langle \xi, w \rangle=    \xi_1 w_1 + \cdots  + \xi_n w_n.
$$
Set
$$
\what{E}_n =  (\R^n \times \R^n \times \C) / \Z^n \times \Z^n.
$$
Now $\what{E}_n = E_1 \otimes \cdots \otimes E_n$, where 
$E_j$ is the pull back of $\what{E}_1$ by the projection $\pi_j:\T^n \times \T^n \to \T \times \T$ onto the $j$-th coordinates.  In \cite{BH2011}, the proof of Theorem 11.3, we showed that 
$\ch(E_j) = 1 +  \beta_j$, where $\beta_j$ is the pull back under $\pi_j$ of the natural generator of $H^2(\T^2;\Z)$.  (Note that the definition of $\what{E}_n $ in \cite{BH2011} has a typo, and it should be defined as above.)   Then the Chern character
$$
 \ch(\what{E}_n) \,\, = \,\, \prod_{j=1}^n \ch(E_j) \,\, = \,\,   \prod_{j=1}^n 1 +  \beta_j,
$$
in particular, it is as non-trivial as it can be. 

Next, note that $\what{E}_n$ restricted to the leaves of $F$ is a flat bundle.  This is because the leaves are products of  $\R$s, where each $\R$ is a leaf of the foliation of the $j$-th $\T^1 \times \T^1$, and the fact that the curvature of the pull back is the pull back of the curvature.  But the connection we will use on $\what{E}_n$ is the tensor product of the pull backs of the  connection on  $\T^1 \times \T^1$ by the $\pi_j$.  So the curvature of our connection on a leaf is the tensor product of the pull backs of the curvatures on one dimensional manifolds, and is identically zero.  

On the open dense subset of $\T^n \times \T^n$ given by $(0,1)^n \times \T^n$,  set 
$$
\nabla(f) \,\, = \,\, df -   (2\pi i \sum_j x_j dw_j)f,
$$
where $f$ is a section of $\what{E}_n$ defined on $(0,1)^n \times \T^n$.  We leave it to the reader to check that in fact this defines a 
connection $\nabla$ on $\what{E}_n$.   $\nabla$  is the tensor product of the pull backs of this connection on $\what{E}_1$ by the $\pi_j$, as promised. 

The leaves of $F$ are given as follows.  Choose $(a,b) \in \R^n \times \R^n$, so that $a_j/b_j \in \R - \Q$ and the  $a_j/b_j $ are not rationally related.  For each $c \in \R^n$, set 
$$
L_c \,\, = \,\, \{a_1 t_1, ... , a_n t_n, b_1  t_1 + c_1, ... , b_n t_n + c_n) \,\, | \, t \in \R^n \}.
$$
The leaves of $F$ are the images of the $L_c$ under the projection $\R^n \times \R^n \to \T^n \times \T^n$, which are denoted $\wtit{L}_c$.   On $(0,1)^n \times \T^n$, we may use $x \in (0,1)^n$ as the coordinates on $\wtit{L}_c$, and when we do, the restriction of $\nabla$ to $\wtit{L}_c$ is  
$$
\nabla^c \,\, = \,\, d   -   2\pi i \sum_j \frac{b_j}{a_j}   x_j dx_j,
$$
which is independent of the transverse coordinate $c$. Note that since the metric on the leaves of $F$ is the usual metric on $\R^n$, the $*$ operator used in the definition of the signature operator is just the usual $*$ operator and is also independent of the transverse coordinate $c$.  So all the elements used in the construction of the twisted leafwise signature operator $D_n$ (which is a generalized Dirac operator) are independent of the transverse coordinate $c$.   Thus, for each $j$, the commutator $[\pa/\pa w_j, D_n] = 0$.  Finally note that the pull back of $\nabla$ to the groupoid $\T^n \times \T^n \times \R^n$ has uniformly bounded coefficients, since we use the usual groupoid coordinates on $\T^n \times \T^n \times \R^n$, which come from $\T^n \times \T^n$.

The operator  $D_n$ is a leafwise operator on the holonomy (homotopy) groupoid for the foliation $F_s$, whose leaves  
are diffeomorphic to the $\wtit{L}_c$, under the target map $r$, with the same structure. Any differential form on $\wtit{L}_c$ twisted by the bundle $\what{E}_n$ can be written as $\omega(t) \otimes s(t)$, where $||s(t)|| = 1$, and
$\nabla^c s(t) = 0$.  An easy computation shows that $D_n$ applied to this is just $D\omega(t) \otimes s(t)$, where $D$ is the usual signature operator on $\R^n$.     As there  are no non-zero $L^2$ harmonic forms on $\R^n$, it follows immediately that $P_0 = 0$ and so is transversely smooth, and that the Chern character of the index bundle $\ch_a([P_0])$ of $D_n$  is zero.

Finally, we show that the spectral projections  $P_{(0,\epsilon)}$ of $D_n$ for the interval  $(0,\epsilon)$  are tranversely smooth.   To this end, we choose the normal bundle $\nu$ of $F$ to be the sub-bundle of $T(\T^n \times \T^n)$ which is just the tangent bundle to the second $\T^n$, that is the bundle spanned by the $\pa/\pa w_j$.  Now the groupoid is diffeomorphic to $\T^n \times \T^n \times \R^n$, and the structure on the leaf $\wL_{(x,w)} = (x,w) \times \R^n$ is given by mapping $(x,w,t)$ to $(at, w+ bt - bx/a)$ in the leaf $L_c$ containing $(x,w)$.  The transverse derivatives of $P_{(0,\epsilon)}$ are obtained by taking vectors  $X \in T(\T^n \times \T^n)_{(x,w)}$, translating them along the leaf $\wL_{(x,w)}$ using the product structure, and then computing the covariant derivatives  $\wtit{\nabla}_X (P_{(0,\epsilon)})$, where  $\wtit{\nabla}$ is the twisted  (using $\what{E}_n$) Levi-Civita connection on $\T^n \times \T^n$.  See \cite{BH2008}, Section 3.  The map from $\wL_{(x,w)}$ to $\wL_{(x  + as ,w + bs)}$, which is parallel translation along the geodesics determined by the vector $\sum_j a_j  \pa /\pa x_j   + b_j \pa /\pa w_j$ in $TF_{(x,w)}$  translated along the leaf $\wL_{(x,w)}$, is actually the identity in our coordinates, so $\wtit{\nabla}_X (P_{(0,\epsilon)}) = 0$ for such $X$.

To see that the $P_{(0,\epsilon)}$ are transversely smooth, note that $ [\pa/\pa w_j,D_n] = 0$ implies that $[\pa/\pa w_j, P_{(0,\epsilon)}] = 0$.   Using the fact that $\wtit{\nabla}$ has uniformly bounded coefficients, we have that $[\wtit{\nabla}_{\pa/\pa w_j}, P_{(0,\epsilon)}]$ is the commutator of the bounded leafwise smoothing operator $P_{(0,\epsilon)}$ with an order zero differential operator with uniformly bounded coefficients, so $[\wtit{\nabla}_{\pa/\pa w_j}, P_{(0,\epsilon)}]$ is a bounded leafwise smoothing operator.  The higher derivatives may be handled the same way, so we have that the spectral projections  $P_{(0,\epsilon)}$  are transversely smooth.

Thus we have that this example satisfies all the conditions of Corollary \ref{limit}, except one.  Namely,  the Novikov-Shubin invariants of the signature operator on $\R^n$ are exactly $n/2$, i.e.\ the Novikov-Shubin invariants of $D_n$ are not greater than half the codimension of the foliation. 

Recall, \cite{HL1999}, Corollary 4, that the Haefliger Connes-Chern character of any leafwise Dirac operator with coefficients in a leafwise flat bundle $E$ is given (up to a constant) by the integral over the fiber of the foliation of the characteristic class $\what{A}(TF)\ch(E)$, where $\what{A}(TF)$ is the $\what{A}$ genus of the tangent bundle of the foliation $F$.   As $TF$ is a trivial bundle,  $\what{A}(TF) = 1$, and the Haefliger Connes-Chern character of 
$D_n$ is given by 
$$
\ch_a(D_n) \,\, = \,\,     \int_F \ch(\what{E}_n).
$$
By Hector et al, \cite{KacimiHector, KacimiHectorSergiescu},   the Haefliger cohomology of $F$ is the same as the basic cohomology , that is the cohomology of the transverse forms which are invariant under the holonomy.  It is not hard to see that this is isomorphic to $H^*(\T^n;\R)$, and we can easily identify the class $\dd \int_F \ch(\what{E}_n)$ under this isomorphism.  In particular it is just $\prod_{j=1}^n  \alpha_j,$  where $\alpha_j$ is the pull back of the natural generator of $H^1(\T^1;\Z)$ under the projection $\T^n \to \T^1$ onto the  $j$\,th coordinate.   Thus $\ch_a(D_n)$ is non-zero, so $\ch_{a}(\Ind_{a}(D)) \neq  \ch_{a}([P_{0}])$ for this example.  Note further that $\ch_a(D_n)$ is non-zero only in the top dimension, so  for $k < n/2$,  $\ch_a^k(D_n) =\ch_a^k([P_0])$,  which they must be by Theorem \ref{main}.

\end{document}